\documentclass[10pt,twoside]{amsart}
\usepackage{lmodern}
\usepackage[T1]{fontenc}
\usepackage{amsmath,amsfonts,amsthm,amssymb,amscd}
\usepackage{latexsym}
\usepackage[pdftex]{graphicx}
\usepackage{enumitem}
\usepackage{bbm}
\usepackage{accents}
\usepackage{euscript}

\usepackage{mathrsfs}

\def\R{ \mathbb{R}}
\def\M{{\mathfrak M}}
\def\m{{\mathfrak m}}
\def\F{ \mathbb{F}}
\def\S{\EuScript{S}}
\newcommand{\lnc}{\mathscr{L}}
\newcommand{\Si}{\mathfrak{S}}

\usepackage{geometry}
\newtheorem{corollary}{Corollary}[section]
\newtheorem{lemma}{Lemma}[section]
\newtheorem{theorem}{Theorem}[section]

\usepackage{hyperref}  
\hypersetup{colorlinks=true, linktoc=all,
    linkcolor=blue,
		citecolor=blue,
		urlcolor=cyan,
		bookmarksopen=true
		}
\usepackage{hyperref}  

\title{On the ternary Estermann problem with almost proportional summands}
\author{Firuz Rakhmonov}
\address{A.Dzhuraev Institute of Mathematics,  National Academy of Sciences of Tajikistan}
\email{rakhmonov.firuz@gmail.com}
\date{}

\begin{document}

\begin{abstract}
For $n \geq 3$, an asymptotic formula is derived for the number of representations of a sufficiently large natural number $N$ in the form $p_1+p_2+m^n=N$, where $p_1$, $p_2$~---~prime numbers, $m$~---~natural number satisfying the conditions
$$
\left|p_k-\mu_kN\right|\le H, \quad k=1,2,\qquad \left|m^n-\mu_3N\right|\le H,\qquad H\ge N^{1-\frac1{n(n-1)}}{\lnc}^{\frac{2^{n+1}}{n-1}+n-1},
$$
for $\mu_1+\mu_2+\mu_3=1, \ \ \mu_i >0. $ \\


Bibliography: 20 titles.
\end{abstract}

\maketitle

\markright{ON THE TERNARY ESTERMANN PROBLEM WITH ALMOST PROPORTIONAL SUMMANDS}

\section{Introduction}
~Estermann \cite{Estermann.eng} proved asymptotic formula for the number of solutions of equation
\begin{equation}\label{formula uravn Estermanna}
 p_1+p_2+m^n=N,
\end{equation}
where $p_1$ and $p_2$ are prime numbers, and $m$ is a natural number, for the case $n=2$. In the works ~\cite{RakhmonovZKh-Matzametki-2003-74-4.eng,RakhmonovZKh-Matzametki-2014-95-3.eng,RAO+RFZ-Iss..Saratov-2016-8.eng} this problem was studied for $n=2,3,4$  under more stringent conditions, specifically when the summands are nearly equal. In these studies, an asymptotic formula was derived for the number of solutions to the Diophantine equation (\ref{formula uravn Estermanna}) under the conditions
$$
\left| p_i-\frac{N}{3}\right|\le H, \quad i=1,2,\qquad \left|m^n-\frac{N}{3}\right|\le H, \qquad H\ge N^{1-\theta(n)}\lnc^{c_n},
$$
corresponding to
\begin{align}\label{formula pannie rezul PrEstcPRS}
\theta(2)=\frac14,\quad c_2=2;\qquad \theta(3)=\frac16,\quad c_3 =3;\qquad\theta(4)=\frac1{12},\quad c_4 =\frac{40}3.
\end{align}

In this paper, a theorem providing an asymptotic formula is obtained for a generalization of Estermann's problem with almost proportional summands for any fixed $n \ge 3$.

{\theorem\label{TeorAsEstermannPPSl} Let  $N$ be a sufficiently large natural number, $n\ge3$ a fixed natural number, and $\rho (N,p)$ the number of solutions to the congruence $x^n\equiv N\pmod p$. Let  $\mu_1$, $\mu_2$, and $\mu_3$ be fixed positive numbers such that $\mu_1+\mu_2+\mu_3=1$, and let $J_n(N,H)$ denote the number of solutions to the Diophantine equation (\ref{formula uravn Estermanna}) under the conditions
$$
\left|p_k-\mu_kN\right|\le H, \quad k=1,2,\qquad \left|m^n-\mu_3N\right|\le H.
$$
Then, for $H\ge N^{1-\frac1{n(n-1)}}\lnc^{\frac{2^{n+1}}{n-1}+n-1},$ the following asymptotic formula holds:
$$
J_{n,r}(N,H)=\frac{3\Si(N)H^2}{ n \mu_3^{1-\frac1n} N^{1-\frac1n} {\lnc}^2 }+O\left(\frac{H^2}{N^{1-\frac1n}\lnc^3}\right),\qquad \Si(N) =\prod_{p}\left(1+\frac{\rho (N,p)}{(p-1)^2}\right)
$$
where the constant in the $O$-term depends on the numbers $\mu_1$, $\mu_2$, $\mu_3$ and $n$.}

From Theorem \ref{TeorAsEstermannPPSl}, when $\mu_1=\mu_2=\mu_3=\frac13$, we obtain an asymptotic formula for the generalization of Estermann's problem with nearly equal summands.
{\corollary \label{TeorAsEstermannPRSl} Let $N$ be a sufficiently large natural number, $n\ge3$ a fixed natural number, and $\rho (N,p)$ the number of solutions to the congruence $x^n\equiv N\pmod p$. Let $J_n(N,H)$ denote the number of solutions to the Diophantine equation (\ref{formula uravn Estermanna}) under the conditions
$$
\left|p_k-\frac{N}3\right|\le H, \quad k=1,2,\qquad \left|m^n-\frac{N}3\right|\le H.
$$
Then, for $H\ge N^{1-\frac1{n(n-1)}}\lnc^{\frac{2^{n+1}}{n-1}+n-1},$ the following asymptotic formula holds:
$$
J_{n,r}(N,H)=\frac{3^{2-\frac1n}\Si(N)H^2}{nN^{1-\frac1n}\lnc^2}+O\left(\frac{H^2}{N^{1-\frac1n}\lnc^3}\right),\qquad \Si(N) =\prod_{p}\left(1+\frac{\rho (N,p)}{(p-1)^2}\right),
$$
where the constant in the $O$-term depends only on $n$.}

It should be noted that the asymptotic formulas obtained earlier in the works~\cite{RakhmonovZKh-Matzametki-2014-95-3.eng,RAO+RFZ-Iss..Saratov-2016-8.eng} for the generalization of Estermann's problem with almost equal summands for $n=3$ and $n=4$, as presented in formula (\ref{formula pannie rezul PrEstcPRS}), are special cases of Corollary \ref{TeorAsEstermannPRSl}.

The proof of Theorem \ref{TeorAsEstermannPPSl}  is carried out using the Hardy-Littlewood-Ramanujan circle method in the form of I.M. Vinogradov's exponential sums, combined with the methods from \cite{RakhmonovZKh-ChebSbornik-2024-25-2.eng,RZKh+RFZ-DNAT-2024-67-3-4.eng,RZKh+RFZ-DNAT-2023-66-9-10.eng,RZKh+RFZ-DNAT-2023-66-11-12.eng,RakhmonovZKh-ChebSbornik-2023-24-3.eng}, where the behavior of short Weyl sums of the form
$$
T(\alpha;x,y)=\sum_{x-y<m\le x}e(\alpha m^n),
$$
 was investigated on magor arcs, and an asymptotic formula for the generalization of Waring's problem with almost proportional summands was proven. The key statements enabling to prove Theorem \ref{TeorAsEstermannPPSl} are:
\begin{itemize}
\item an asymptotic formula for short Weyl trigonometric sums of the form $T(\alpha; x, y)$ in a small neighborhood of the centers of the major arcs (Corollary~\ref{Sledst1 lemmi ob poved kor trig summi Weyl} to Lemma~\ref{Lemma ob poved kor trig summi Weyl});
\item a nontrivial estimate for the sums $T(\alpha; x, y)$ on the major arcs outside small neighborhoods of their centers (Corollary~\ref{Sledst2 lemmi ob poved kor trig summi Weyl} to Lemma~\ref{Lemma ob poved kor trig summi Weyl});
\item a nontrivial estimate for the sums $T(\alpha; x, y)$ on the minor arcs (Theorem~\ref{Teor otsenka kor summ Weyl n>=3}).
\end{itemize}

When solving a number of additive problems with almost proportional summands, such as Waring's problem and the generalization of Estermann's problem, issues related to the behavior of short Weyl exponential sums $T(\alpha;x,y)$ arise both on magor and minor arcs. For arbitrary $n\ge3 $, these sums on large arcs were studied in \cite{RakhmonovZKh-ChebSbornik-2024-25-2.eng,RZKh+RFZ-DNAT-2023-66-11-12.eng}. Non-trivial estimates for the exponential sums $T(\alpha;x,y)$ on minor arcs for $q\gg y^\varepsilon$ were obtained in \cite{RZKh+AAZ+NNN-DANRT-2018-7-8.eng}.

In this work, using the methods from \cite{RakhmonovFZ-ChebSbornik-2011-12-1.eng,RakhmonovFZ-VestnikMGU-2011-3.eng,RakhmonovZKh+RFZ-TrMIRAN-2016-296.eng,RakhmonovZKh+FZ-ChebSbornik-2019-20-4.eng},  a non-trivial estimate is obtained under the condition
$$
(\ln y)^{(n-1)^2}\ll q\ll y^n(\ln y)^{-(n-1)^2}.
$$
{\theorem\label{Teor otsenka kor summ Weyl n>=3} Let $x\ge x_0>0$, $\ln x<y\le x(\ln x)^{-1}$,  and $\alpha$ be a real number,
$$
\left|\alpha-\frac{a}{q}\right|\le \frac{1}{q^2}, \qquad (a,q)=1,
$$
then for  $n\ge3$ the following estimate holds:
$$
|T(\alpha;x,y)|\ll y\left(\frac1q+\frac1y+\frac{q}{y^n}\right)^{2^{-n}}(\ln qy)^{(n-1)^22^{-n}},
$$
where the implied constant in Vinogradov's symbol $\ll$ depends only on $n$.}

\emph{Notation.} Let $\varepsilon $ be any positive number not exceeding $0.00001$,  ${\mathscr L}=\ln N$,
\begin{align*}
&S(a,q)=\sum_{k=1}^qe\left(\frac{ak^n}q\right),
&\gamma(\lambda;x,y)=\int_{-0,5}^{0,5}e\left(\lambda\left(x-\frac y2+yt\right)^n\right)dt.
\end{align*}

\section{Auxiliary Statements}
 {\lemma \label{Lemma virazh st T(f(u);x,y) konech razn} {\rm \cite{RakhmonovZKh-ChebSbornik-2024-25-2.eng}}. Let $f(u)$ be a polynomial of degree $n$, $x$ and $y$ be positive integers with $y<x$. Then, for $j=1,\ldots,n-1$ we have the inequality
\begin{align*}
&\left|\sum_{x-y<u\le x}e(f(u))\right|^{2^{j}}\leq(2y)^{2^j-j-1}\sum_{|h_1|<y}\ldots\sum_{|h_j|<y}\left|\sum_{u\in I_j(x,y;h_1,\ldots,h_j)}e(\Delta_j(f(u); h_1,\ldots,h_j))\right|,
\end{align*}
where the intervals $I_j(x,y;h_1,\ldots,h_j)$ are defined by the relations:
\begin{align*}
&I_1(x,y;h_1)=(x-y,x]\cap (x-y-h_1,x-h_1],\\
&I_j(x,y;h_1,\ldots, h_{j})=I_{j-1}(x,y;h_1,\ldots, h_{j-1})\cap I_{j-1}(x-h_j,y;h_1,\ldots, h_{j-1}),
\end{align*}
i.e., $I_{j-1}(x-h_j,y;h_1,\ldots, h_{j-1})$ is obtained from $I_{j-1}(x,y;h_1,\ldots, h_{j-1})$ by shifting all the intervals that constitute it as an intersection, by $-h_j$, .}

{\lemma \label{Lemma razn oper} {\rm \cite{RakhmonovZKh-ChebSbornik-2024-25-2.eng}}. Let  $\Delta_j$ denote the $j$-th application of the difference operator, so for any real-valued function $f(u)$ we have:
\begin{align*}
&\Delta_1(f(u);h)=f(u+h)-f(u),\nonumber\\
&\Delta_{j+1}(f(u);h_1,\ldots,h_{j+1})=\Delta_1(\Delta_k(f(u);h_1,\ldots,t_j);h_{j+1}).
\end{align*}
Then, for $j=1,\ldots,n-1$, the following relation holds:
$$
\Delta_j(u^n;h_1,\ldots,h_j)=h_1\ldots h_jp_j(u;h_1,\ldots,h_j),
$$
where $p_j=g_j(u;h_1,\ldots,h_j)$ is a form of degree $n-j$ with integer coefficients, having degree $n-j$ with respect to $u$, and its leading coefficient is $n(n-1)\ldots(n-j+1)$, i.e.,
$$
p_j(u;h_1,\ldots,h_j)=\frac{n!}{(n-j)!}u^{n-j}+\ldots.
$$}
{\lemma \label{Lemma ob otsenke summ st obsh funk delit} {\rm \cite{Mardjanishvili-garv}.} For $x\ge 1$, $r\ge2$  and $k\ge1$, the following estimate holds:
\begin{align*}
\sum_{n\le x}\tau_r^k(n)\ll\frac{xr^k}{(r!)^{\frac{r^k-1}{r-1}}}\left(\ln x+r^k-1\right)^{r^k-1}.
\end{align*}}
{\lemma \label{Lemma ob otsenke summ obr RDBTsCh} {\rm \cite{Karatsuba-OATCh.eng}.}  Let $\alpha$ be a real number,
$$
\left|\alpha-\frac{a}{q}\right|\le \frac{1}{q^2}, \qquad (a,q)=1,
$$
with $x\ge1$, $y>0$  and any $\beta$, then the following estimate holds:
$$
\sum_{n\le x}\min\left(y,\ \frac1{\|\alpha n+\beta\|}\right)\le6\left(\frac{x}q+1\right)(y+q\ln q).
$$}
{\lemma \label{Lemma KLTrSummPr-MalOkrSenBolDug} {\rm \cite{SAA-DNAT-2021-64-11-12.eng}}. Let $x\ge x_0$, $A$ and $b$ be arbitrary fixed positive numbers, $1\le q\le\lnc_x^b$,
$$
\alpha=\frac{a}{q}+\lambda,\quad (a,q)=1.
$$
Then, under the conditions  $|\lambda|\le x\left(2\pi y^2\right)^{-1}$ and  $y\ge x^{\frac{\scriptstyle5}{\scriptstyle8}}\lnc_x^{1,5A+0,25b+18}$, the following equality holds:
$$
S_1(\alpha;x,y)=\frac{\mu(q)}{\varphi(q)}\frac{\sin\pi\lambda y}{\pi\lambda}e\left(\lambda\left(x-\frac{y}{2}\right)\right) +O\left(y\lnc_x^{-A}\right).
$$}
{\lemma \label{Lemma S(alpha;Nk,H=S1(alpha;Nk+H,2H)} Let $\mu_k$ be a fixed real number,  $0<\mu_k<1$, $N$ be a sufficiently large natural number, $N_k=\mu_kN+H$, $k=1,\ 2$, $N^\frac12\le H\le N^{1-\frac1{30}}$,
$$
\S(\alpha;N_k,2H)=\sum_{N_k-2H<p\le N_k}e(\alpha p),\qquad S_1(\alpha;x,y)=\sum_{x-y<n\le x}\Lambda(n)e(\alpha n).
$$
Then the following relation holds:
$$
\S(\alpha;N_k,2H)=\frac{S_1(\alpha;N_k,2H)}{\ln(\mu_kN)}+O\left(\frac{H^2}{N\ln(\mu_kN)}\right).
$$}

{\sc Proof:} Expressing the inequality $N_k-2H<p\le N_k$   in the form $\mu_kN-H<p\le\mu_kN+H$, taking the logarithm, transforming this inequality and using the formula:
$$
\ln(\mu_kN\pm H)=\ln(\mu_kN)+\ln\left(1\pm\frac{H}{\mu_kN}\right)=\ln(\mu_kN)+O\left(\frac{H}{N}\right),
$$
we obtain that for  $\mu_kN-H<p\le\mu_kN+H$, the following relation holds:
\begin{align*}
& \ln p=\ln(\mu_kN)+O\left(\frac HN\right). 
\end{align*}
Using this formula, the sum $\S(\alpha;N_k,2H)$ can be expressed in terms of the sum $S_1(\alpha;x,y)$. We have:
\begin{align}
\S(\alpha;N_k,2H)&=\sum_{\mu N-H<p\le\mu_kN+H}\left(\frac{\ln p}{\ln(\mu_kN)}+O\left(\frac{H}{N\ln(\mu_kN)}\right)\right)e(\alpha p)=\nonumber\\
&=\frac{S_1(\alpha;\mu_kN+H,2H)}{\ln(\mu_kN)}-R_1+O\left(\frac{H^2}{N\ln(\mu_kN)}\right). \label{Formula-S1(alpha;muiN,H)=S1(alpha;muiN+H,2H)...}
\end{align}
Estimating the sum  $R_1$ trivially by the number of terms and using the formula
$$(1\pm u)^{\mu_k}=1\pm\mu_ku+O(u^2), \ |u|<0,5,$$ we get:
\begin{align*}
R_1&=\sum_{\substack{\mu_kN-H<p^k\le\mu_kN+H\\ k\ge2}}\frac{\ln p}{\ln(\mu_kN)}\ e(\alpha p^{k})\ll\ln(\mu_kN)\left(\left(\mu_kN+H\right)^\frac12-\left(\mu_kN-H\right)^\frac12+1\right)=\\
&=(\mu_kN)^\frac12\ln(\mu_kN)\left(\left(1+\frac{H}{\mu_kN}\right)^\frac12-\left(1-\frac{H}{\mu_kN}\right)^\frac12\right)+\ln(\mu_kN)
\ll\frac{H}{(\mu_kN)^\frac12}\ln(\mu_kN).
\end{align*}
Substituting this estimate into the right-hand side of (\ref{Formula-S1(alpha;muiN,H)=S1(alpha;muiN+H,2H)...}), we obtain the statement of the lemma.
{\lemma \label{Lemma ob poved kor trig summi Weyl} {\rm \cite{RakhmonovZKh-ChebSbornik-2024-25-2.eng}}. Let $\tau\geq 2n(n-1)x^{n-2}y$, then for  $\{n|\lambda|x^{n-1}\}\le \frac{1}{2q}$ the following formula holds:
\begin{equation*}
T(\alpha;x,y)=\frac{S(a,q)}{q}T(\lambda;x,y)+O\left(q^{\frac{1}{2}+\varepsilon}\right),
\end{equation*}
and for $\{n|\lambda|x^{n-1}\}>\frac{1}{2q}$ we have the estimate
\begin{align*}
|T(\alpha,x,y)|&\ll q^{1-\frac 1n}\ln q +\min(yq^{-\frac 1n}, \lambda^{-\frac12} x^{1-\frac n2}q^{-\frac1n}).
\end{align*}}

{\corollary \label{Sledst1 lemmi ob poved kor trig summi Weyl} Let $\tau\geq 2n(n-1)x^{n-2}y$,  $|\lambda|\le \frac{1}{2nq x^{n-1}}$, then the relation
\begin{align*}
T(\alpha;x,y)=\frac{y}{q}S(a,q)\gamma(\lambda;x,y)+O(q^{\frac 12+\varepsilon})
\end{align*}}
holds.
{\corollary \label{Sledst2 lemmi ob poved kor trig summi Weyl} Let $\tau\geq2n(n-1)x^{n-2}y$, $\frac{1}{2nq x^{n-1}}<|\lambda|\le\frac1{q\tau}$, then the estimate
$$
T(\alpha;x,y)\ll q^{1-\frac1n}\ln q +\min\left(yq^{-\frac1n}, x^\frac12q^{\frac12-\frac1n}\right)
$$} holds.
{\lemma \label{Lemma otsenka polnoy trigsumm}{\rm \cite{AGI+KAA+ChVN-TKTS-1987-Nauka.eng}}. Let $(a,q)=1$, $q$ be a natural number. Then we have:
$$
S(a,q)=\sum_{k=1}^{q}e\left(\frac{ak^n}{q}\right)\ll q^{1-\frac1n},
$$
where the constant in Vinogradov's symbol depends on $n$.}
{\lemma \label{Lemma Huxley pi(x)-pi(x-y)}{\rm \cite{Huxley1972-eng}.} Let $y\ge x^{\frac{7}{12}+\varepsilon}$, then the following asymptotic formula holds:
$$
\pi(x)-\pi(x-y)=\frac{y}{\ln x} +O\left(\frac{y}{\ln^2 x}\right).
$$}
{\lemma \label{Lemma otsenka trig integ po pervoy proizv}{\rm \cite{AGI+KAA+ChVN-TKTS-1987-Nauka.eng}}. Let the real function $f(u)$ and the monotonic function $g(u)$ satisfy the conditions: $f'(u)$ is monotonic, $|f'(u)|\ge m_1>0$ and $|g(u)|\le M$. Then the following estimate holds:
$$
\int_a^bg(u)e(f(u))du\ll \frac{M}{m_1}.
$$}
{\lemma \label{Lemma otsenka trig int po n-oy proizvodn}{\rm \cite{AGI+KAA+ChVN-TKTS-1987-Nauka.eng}}.  Let $f(u)$ be a real function for $a\le u\le b$, with $n$-th order derivative $f^{(n)}(u)$ (where $n$>1), and for some $A>0$ the inequality $A\le |f^{(n)}(u)|$ holds. Then the following estimate holds:
$$
\int\limits_a^be(f(u))du\le \min (b-a,6nA^{-\frac 1n}).
$$}

\section{Proof of Theorem \ref{Teor otsenka kor summ Weyl n>=3}}
Using Lemma \ref{Lemma virazh st T(f(u);x,y) konech razn} for $j=n-1$, and then applying Lemma \ref{Lemma razn oper}, we obtain
\begin{align*}
|T(\alpha;x,y)|^{2^{n-1}}&\leq(2y)^{2^{n-1}-n}\sum_{|h_1|<y}\ldots\sum_{|h_{n-1}|<y}\left|\sum_{u\in I_{n-1}(x,y;h_1,\ldots,h_{n-1})}e(\alpha n!h_1h_2\ldots h_{n-1}u)\right|.
\end{align*}
In the last sum over $u$, the number of terms for which the relation $h_1\cdot\cdot\cdot h_{n-1}=0$ holds, does not exceed $(n-1)y(2y)^{n-2}$. Therefore, we have
\begin{align}\label{formula2}
|T(\alpha;&x,y)|^{2^{n-1}}\le(2y)^{2^{n-1}-n}\left(2^{n-1}\tilde{T}(\alpha;x,y)+(n-1)y(2y)^{n-2}\right),\\
&\tilde{T}(\alpha;x,y)=\sum_{1\le h_1<y}\ldots\sum_{1\le h_{n-1}<y}\left|\sum_{u\in I_{n-1}(x,y;h_1,\ldots,h_{n-1})}e(\alpha n!h_1h_2\ldots h_{n-1}u)\right|\le\nonumber\\
&\le\sum_{h=1}^{y^{n-1}}\tau_{n-1}(h)\left|\sum_{u\in I_{n-1}(x,y;h_1,\ldots,h_{n-1})}e(\alpha n!hu)\right|\le\sum_{h=1}^{y^{n-1}}\tau_{n-1}(h)\min\left(y,\frac1{2\|\alpha n!h\|}\right),\nonumber
\end{align}
where $\tau_n(h)$ is the number of solutions to the Diophantine equation $h_1\ldots h_n=h$.
Applying the Cauchy inequality to the last sum, followed by Lemmas \ref{Lemma ob otsenke summ st obsh funk delit} and \ref{Lemma ob otsenke summ obr RDBTsCh}, we find
\begin{align*}
\tilde{T}^2(\alpha;x,y)&\le\sum_{h=1}^{y^{n-1}}\tau_{n-1}^2(h)\cdot y\sum_{h=1}^{y^{n-1}}\min\left(y,\frac{1}{2\|\alpha n!h\|}\right)\le\\
&\le\frac{(n-1)^2y^n}{((n-1)!)^2}(\ln y-n^2-2n)^{n^2-2n}\left(\frac{n!y^{n-1}}q+1\right)(y+q\ln q)\ll\\
&\ll y^{2n}\left(\frac1q+\frac1y+\frac{q\ln q}{y^n}\right)(\ln y)^{n^2-2n}.
\end{align*}
Substituting this estimate into (\ref{formula2}), we obtain
\begin{align*}
&|T(\alpha;x,y)|^{2^n}\ll y^{2^n-2n}\left(\tilde{T}^2(\alpha;x,y)+y^{2n-2}\right)\ll \\
&\ll y^{2^n-2n}\left(y^{2n}\left(\frac1q+\frac1y+\frac{q\ln q}{y^n}\right)(\ln y)^{n^2-2n}+y^{2n-2}\right)\ll \\
&\ll y^{2^n}\left(\frac1q+\frac1y+\frac{q}{y^n}\right)(\ln qy)^{(n-1)^2}.
\end{align*}
Taking the $2^n$-th root yields the assertion of the theorem.

\section{Proof of Theorem \ref{TeorAsEstermannPPSl}}
Without loss of generality, we assume that
\begin{equation}\label{formula opred teta i omega}
H=N^{1-\theta(n)}\lnc^{\omega(n)},\qquad \theta(n)=\frac1{(n-1)n},\qquad \omega(n)=\frac{2^{n+1}}{n-1}+n-1.
\end{equation}
Using the notations
\begin{align}
\S(\alpha;&N_k,2H)=\sum_{N_k-2H<p\le N_k}e(\alpha p),\qquad N_k=\mu_kN+H,\qquad k=1,\ 2;\nonumber  \\
&N_3=(\mu_3N+H)^\frac1n=\mu_3^\frac1nN^\frac1n\left(1+\frac{H}{\mu_3N}\right)^\frac1n=\mu_3^\frac1nN^\frac1n\left(1+O\left(\frac HN\right)\right), \label{formula poryadok N3}\\
&H_3=(\mu_3N+H)^\frac1n-(\mu_3N-H)^\frac1n=\frac{2H}{n\mu_3^{1-\frac1n}N^{1-\frac1n}}\left(1+O\left(\frac{H^2}{N^2}\right)\right), \label{formula poryadok H3}\\
&\tau=2(n-1)nN_3^{n-2}H_3=\frac{4(n-1)H}{n\mu_3^\frac1nN^\frac1n}\left(1+O\left(\frac HN\right)\right), \qquad \varkappa \tau=1, \label{formula poryadok tau}
\end{align}
$J_n(N,H)$ -- denotes the number of solutions to the Diophantine equation
$$
p_1+p_2+m^n=N,
$$
with the prime numbers $p_1$, $p_2$ and the natural numbers $m$ under the conditions
$$
|p_k-\mu_kN|\le H, \quad k=1,2,\quad |m^n-\mu_3N|\le H, \quad \mu_1+\mu_2+\mu_3=1, \quad \mu_i>0,
$$
which can be expressed as
\begin{align}
J_n(N,H)&=\int_{-\varkappa}^{1-\varkappa}e(-\alpha N)\sum_{|p_1-\mu_1 N|\le H}e(\alpha p_1)\sum_{|p_2-\mu_2 N|\le H}e(\alpha p_2)\sum_{|m^n-\mu_3N|\le H}e(\alpha m^n)=\nonumber\\
&=\int_{-\varkappa}^{1-\varkappa}e(-\alpha N)\left(\S(\alpha;N_1,2H)+\theta_1\right)\left(\S(\alpha;N_2,2H)+\theta_2\right) \left(T(\alpha;N_3,H_3)+\theta_3\right)d\alpha,\label{formula Jn(N,H)-1}
\end{align}
where $|\theta_k|$ equals $1$ if the lower bounds of the exponential sums $\S(\alpha;N_1,2H)$, $\S(\alpha;N_2,2H)$, $T(\alpha;N_3,H_3)$, i.e., the numbers $N_1-2H$, $N_2-2H$, $N_3-H_3$, are integers, and $0$ - otherwise. Multiplying the brackets in the integrand in (\ref{formula Jn(N,H)-1}), we obtain
\begin{align}
J_n(N,H)=\int_{-\varkappa}^{1-\varkappa}e(-\alpha N)&\S(\alpha;N_1,2H)\S(\alpha;N_2,2H)T(\alpha;N_3,H_3)d\alpha+\R_1,\label{formula Jn(N,H)-2}\\
\R_1=\int_{-\varkappa}^{1-\varkappa}e(-\alpha N)&\left(
\theta_3\S(\alpha;N_1,2H)\S(\alpha;N_2,2H)+ \theta_1\theta_2T(\alpha;N_3,H_3)+\right.\nonumber\\
&\left.+\theta_2\S(\alpha;N_1,2H)T(\alpha;N_3,H_3)+\theta_1\theta_3\S(\alpha;N_2,2H)+\right. \nonumber\\
&\left.+\theta_1\S(\alpha;N_2,2H)T(\alpha;N_3,H_3)+\theta_2\theta_3\S(\alpha;N_1,2H)\right)d\alpha.\nonumber
\end{align}
 In $\R_1$, moving on to estimates, using the Cauchy inequality and the following relations:
\begin{align*}
&\int_0^1|\S(\alpha;N_k,2H)|^2d\alpha=\pi(N_k)-\pi(N_k-2H)\le 2H, \qquad k=1,\ 2;\\
&\int_0^1|T(\alpha;N_3,H_3)|^2d\alpha=[N_3]-[N_3-H_3]\le H_3+1\ll\frac{H}{N^{1-\frac1n}},
\end{align*}
where the second relation utilizes the relation (\ref{formula poryadok H3}), we have
\begin{align*}
\R_1&\le\left(\int_0^1|\S(\alpha;N_1,2H)|^2d\alpha\int_0^1|\S(\alpha;N_2,2H)|^2d\alpha\right)^\frac12+ \left(\int_0^1|T(\alpha;N_3,H_3)|^2d\alpha\right)^\frac12+\\
&+\left(\int_0^1|\S(\alpha;N_1,2H)|^2d\alpha\int_0^1|T(\alpha;N_3,H_3)|^2d\alpha\right)^\frac12+ \left(\int_0^1|\S(\alpha;N_2,2H)|^2d\alpha\right)^\frac12+\\
&+\left(\int_0^1|\S(\alpha;N_2,2H)|^2d\alpha\int_0^1|T(\alpha;N_3,H_3)|^2d\alpha\right)^\frac12+ \left(\int_0^1|\S(\alpha;N_1,2H)|^2d\alpha\right)^\frac12\ll\\
&\ll H+ \left(\frac{H}{N^{1-\frac1n}}\right)^\frac12+\left(\frac{H^2}{N^{1-\frac1n}}\right)^\frac12+H^\frac12\ll H\ll\frac{H^2}{N^{1-\frac1n}\lnc^3}.
\end{align*}
From this and from (\ref{formula Jn(N,H)-2}), we obtain
$$
J_n(N,H)=\int_{-\varkappa}^{1-\varkappa}e(-\alpha N)\S(\alpha;N_1,2H)\S(\alpha;N_2,2H)T(\alpha;N_3,H_3)d\alpha+O\left(\frac{H^2}{N^{1-\frac1n}\lnc^3}\right).
$$
According to Dirichlet's theorem on the approximation of real numbers by rational numbers, every $\alpha $ in the interval $[-\varkappa,1-\varkappa]$ can be expressed in the form
\begin{equation}
\label{predst alpha=a/q+lambda,...}
 \alpha=\frac aq+\lambda,\qquad (a,q)=1,\qquad 1\le q\le\tau, \qquad |\lambda|\le\frac1{q\tau}.
\end{equation}
In this representation $0\le a\le q-1$, with $a=0$ only when $q=1$. Let $\M$ denote those $\alpha $ for which in the representation (\ref{predst alpha=a/q+lambda,...}) we have
$$
q\le\lnc^{\eta(n)},\qquad \eta(n)=2^{n+1}+(n-1)^2,
$$
and let $\m$ denote the remaining $\alpha$. The set $\M$ consists of disjoint intervals.  We will partition the set $\M$ into subsets $\M_1$ and $\M_2$:
\begin{align*}
&\M_1=\left\{\alpha:\quad \alpha\in\M,\quad \left|\alpha-\frac aq\right|\le\frac{\lnc^2}H\right\} ,\\
&\M_2=\left\{\alpha:\quad \alpha\in\M,\quad \frac{\lnc^2}H<\left|\alpha-\frac aq\right|\le\frac1{q\tau}\right\} .
\end{align*}
 Let $J(\M_1)$, $J(\M_2)$ and $J(\m)$ denote the integrals over the sets $\M_1$, $\M_2$ and $\m$, respectively. We have
$$
J_n(N,H)=J(\M_1)+J(\M_2)+J(\m)+O\left(\frac{H^2}{N^{1-\frac1n}\lnc^3}\right).
$$
In the last formula, the first term, $J(\M_1)$, provides the main term of the asymptotic formula for $J(N,H)$, while $J(\M_2)$ and $J(\m)$ contribute to its remainder term.

 \subsection{Calculation of the Integral $J(\M_1)$.} By definition, the integral $J(\M_1)$ is given by:
\begin{align}\label{predst I(M_1)= sum sumI(a,q)}
J(\M_1)&=\int_{\M_1}\F(\alpha)e(-\alpha N)d\alpha=\sum_{q\le\lnc^{\eta(n)}}\sum_{\substack{a=0\\ (a,q)=1}}^{q-1}I(a,q), \\
& I(a,q)=e\left(-\frac{aN}q\right)\int_{|\lambda|\le\lnc^2H^{-1}}\F\left(\frac aq+\lambda;N,H\right)e(-\lambda N)d\lambda ,\label{I(a,q)=intF(alpha;N,H)} \\
& \F(\alpha)=\F(\alpha;N,H)=\S(\alpha;N_1,2H)\S(\alpha;N_2,2H)T(\alpha;N_3,H_3). \nonumber
\end{align}
To derive the asymptotic formula for the function $\F(\alpha;N,H)$, we first determine the asymptotic behavior of the sum $S_1(\alpha;N_k,2H)$, $k=1;\ 2$. We set:
$$
x=\mu_kN+H,\qquad y=2H,\qquad A=3+\frac{n-1}n\eta(n),\qquad b=\eta(n),
$$
and apply Lemma \ref{Lemma KLTrSummPr-MalOkrSenBolDug} to these sums. Utilizing the relations
$$
1,5A+0,25\eta(n)+18=\frac{45}2+\left(\frac54-\frac1n\right)\eta(n),\qquad 1-\theta(n)\ge1-\theta(3)=\frac56<\frac58,
$$
we show that the following conditions of Lemma \ref{Lemma KLTrSummPr-MalOkrSenBolDug} are satisfied, namely:
\begin{align*}
&2H=2N^{1-\theta(n)}\lnc^{\omega(n)}\ge2N_k^{1-\theta(n)}\lnc^{\omega(n)}\ge N_k^{\frac{\scriptstyle5}{\scriptstyle8}}\left(\ln N_k\right)^{\frac{45}2+\left(\frac54-\frac1n\right)\eta(n)}, \\
&\frac{\lnc^2}{H}=\frac{\mu_k N+2H}{2\pi(2H)^2}\cdot\frac{8\pi H\lnc^2}{\mu_k N+2H}\le\frac{\mu_k N+2H}{2\pi(2H)^2}.
\end{align*}
Accordingly, we obtain:
\begin{align*}
S_1(\alpha;N_k,2H )&=\frac{\mu(q)}{\varphi(q)}\frac{\sin2\pi\lambda H }{\pi\lambda}e\left(\lambda\mu_kN\right) +O\left(\frac{H}{\left(\ln(\mu_kN)\right)^A}\right).
\end{align*}
Using Lemma \ref{Lemma S(alpha;Nk,H=S1(alpha;Nk+H,2H)}, the sum $\S(\alpha;N_k,2H)$, $k=1,\ 2$ can be expressed through the sum $S_1(\alpha;N_k,2H)$, and then, using the last formula, we have
\begin{align*}
\S(\alpha;\mu_kN,H)&=\frac{S_1(\alpha;\mu_kN+H,2H)}{\ln(\mu_kN)}+O\left(\frac{H^2}N\right)= \frac{\mu(q)}{\varphi(q)}\frac{\sin2\pi\lambda H}{\pi\lambda}\frac{e(\lambda\mu_kN)}{\ln(\mu_kN)}+O\left(\frac{H}{\lnc^{A+1}}\right).
\end{align*}
From this, and the relation $\dfrac{\sin2\pi\lambda H}{\pi\lambda}\ll H$, we find
\begin{align}
\S(\alpha;\mu_1N,H)\S(\alpha;\mu_2N,H)=&\frac{\mu^2(q)}{\varphi^2(q)}\frac{\sin^22\pi\lambda H }{\pi^2\lambda^2}\frac{e(\lambda(\mu_1+\mu_2)N)}{\ln(\mu_1N)\ln(\mu_2N)}+O\left(\frac{H^2}{\varphi(q)\lnc^{A+2}}\right). \label{formula asimp S32(alpha;N1,2H1)}
\end{align}
Now, using Corollary \ref{Sledst1 lemmi ob poved kor trig summi Weyl} of Lemma \ref{Lemma ob poved kor trig summi Weyl}, we derive the asymptotic behavior of the sum $T(\alpha;N_3,H_3)$, $\alpha\in\M_1$. From the inequality $4nN^{1-\frac1n}\lnc^{\eta(n)+2}\le H=N^{1-\theta(n)}\lnc^{\omega(n)}$ , it follows that the conditions of this lemma are satisfied, namely
\begin{align*}
\frac{\lnc^2}{H}=&\frac1{2nqN_3^{n-1}}\cdot\frac{2nqN_3^{n-1}\lnc^2}{H}=\frac1{2nqN_3^{n-1}}\cdot
\frac{2nq\mu_3^\frac{n-1}nN^{1-\frac1n}\lnc^2}H\left(1+O\left(\frac HN\right)\right)\le\\
&\le\frac1{2nqN_3^{n-1}}\cdot\frac{4nN^{1-\frac1n}\lnc^{\eta(n)+2}}H\le\frac1{2nqN_3^{n-1}}.
\end{align*}
Therefore, according to Corollary \ref{Sledst1 lemmi ob poved kor trig summi Weyl} of Lemma \ref{Lemma ob poved kor trig summi Weyl}, we obtain
\begin{equation}\label{formula asimpT(alpha;N3,H3)-1}
T(\alpha;N_3,H_3)=\frac{H_3}{q}S(a,q)\gamma(\lambda;N_3,H_3)+O(q^{\frac12+\varepsilon}).
\end{equation}
Now, let us derive the asymptotic formula for $\gamma(\lambda;N_3,H_3)$. Using the relation $N_3^n=\mu_3N+H$, and the formula
$$
N_3^{n-i}H_3^i=\frac{2^iN^{1-i}H^i}{n\mu_3^{i-1}}\left(1+O\left(\frac HN\right)\right),
$$
which follows from relations (\ref{formula poryadok N3})  and (\ref{formula poryadok H3}), we have
\begin{align*}
&\left(N_3+H_3\left(t-\frac12\right)\right)^n=N_3^n+nN_3^{n-1}H_3\left(t-\frac12\right) +\sum_{i=2}^nC_n^iN_3^{n-i}H_3^i\left(t-\frac12\right)^i=\\
&=\mu_3N+H+(2Ht-H)\left(1+O\left(\frac HN\right)\right)+\sum_{i=2}^nC_n^i\frac{2^iN^{1-i}H^i}{n\mu_3^{i-1}} \left(t-\frac12\right)^i\left(1+O\left(\frac HN\right)\right)=\\
&=\mu_3N+2Ht+\R_2,\qquad \R_2\ll\frac{H^2}N.
\end{align*}
Thus, noting that $e(\lambda\R_2)-1\ll|\lambda|\R_2$ and $|\lambda|\ll\lnc^2H^{-1}$, we find
\begin{align*}
&\gamma(\lambda;N_3,H_3)=\int_{-0,5}^{0,5}e\left(\lambda\left(\mu_3N+2Ht+\R_2\right)\right)dt
=e\left(\mu_3N\lambda\right)\frac{\sin\left(2\pi H\lambda\right)}{2\pi H\lambda}+O\left(\frac{H\lnc^2}N\right).
\end{align*}
Substituting the right-hand side of this formula into (\ref{formula asimpT(alpha;N3,H3)-1}), and then using the estimate for the remainder term $\R_3$ obtained by applying relation (\ref{formula poryadok H3}), the estimate $S(a,q)\ll q^{1-\frac1n}$ (Lemma \ref{Lemma otsenka polnoy trigsumm}), and the condition $q\le\lnc^{\eta(n)}$, as well as relation (\ref{formula opred teta i omega}), we obtain
\begin{align}
&T(\alpha;N_3,H_3)=\frac{S(a,q)}{q}\frac{H_3\sin\left(2\pi H\lambda\right)}{2\pi H\lambda}(e(\mu_3N\lambda)+\R_3,\label{formula asimpT(alpha;N3,H3)-2}\\
\R_3&\ll\frac{H^2\lnc^2}{q^\frac1nN^{2-\frac1n}}+q^{\frac12+\varepsilon} =\frac{H^2\lnc^2}{q^\frac1nN^{2-\frac1n}}\left(1+\frac{N^{2-\frac1n}\lnc^{-2}q^{\frac12-\frac12n+\varepsilon}}{H^2}\right)=\nonumber\\
&=\frac{H^2\lnc^2}{q^\frac1nN^{2-\frac1n}}\left(1+N^{-\frac{n-3}{n(n-1)}}\lnc^{-2-2\omega(n)}q^{\frac12-\frac1n +\varepsilon}\right)\le\nonumber\\
&\le\frac{H^2\lnc^2}{q^\frac1nN^{2-\frac1n}}\left(1+N^{-\frac{n-3}{n(n-1)}}\lnc^{\nu(n)}\right), \quad \nu(n)=-2-2\omega(n)+\left(\frac12-\frac1n +\varepsilon\right)\eta(n),  \nonumber
\end{align}
From here and taking into account that $\nu(3)=-\frac{56}3+20\varepsilon$, we find
\begin{align*}
\R_3\ll\frac{H^2\lnc^2}{q^\frac1nN^{2-\frac1n}}.
\end{align*}
By term-by-term multiplication of formulas (\ref{formula asimp S32(alpha;N1,2H1)}) and (\ref{formula asimpT(alpha;N3,H3)-2}), and then estimating the remainder term of this product, denoted by $\R_4$, using the inequality $|\sin(2\pi H\lambda)|\ll |2\pi H\lambda|$, the estimate $S(a,q)\ll q^{1-\frac1n}$ (Lemma \ref{Lemma otsenka polnoy trigsumm}), relation (\ref{formula poryadok H3}), and condition (\ref{formula opred teta i omega}), we obtain:
\begin{align*}
\F(\alpha)&=\frac{H_3}{2H\ln(\mu_1N)\ln(\mu_2N)}\cdot\frac{\mu^2(q)S(a,q)}{q\varphi^2(q)}\cdot\frac{\sin^32\pi\lambda H}{\pi^3\lambda^3}+\R_4,\\
\R_4&\ll\frac1{\varphi^2(q)}\frac{\sin^22\pi\lambda H }{\pi^2\lambda^2\lnc^2}\R_3+ \frac{|S(a,q)|}{q}\frac{H_3|\sin(2\pi H\lambda)|}{|2\pi H\lambda|}\frac{H^2}{\varphi(q)\lnc^{A+2}}+\frac{H^2}{\varphi(q)\lnc^{A+2}}\R_3\ll\\
&\ll\frac{H^4}{\varphi^2(q)q^\frac1nN^{2-\frac1n}}+\frac{H^3}{q^\frac1n\varphi(q)N^{1-\frac1n}\lnc^{A+2}}+ \frac{H^4}{q^\frac1n\varphi(q)N^{2-\frac1n}\lnc^A}=\\
&=\frac{H^3}{q^\frac1n\varphi(q)N^{1-\frac1n}\lnc^{A+2}}\left(1+\frac{H\lnc^{A+2}}{\varphi(q)N}+\frac{H\lnc^2}{N}\right) \ll\frac{H^3}{q^\frac1n\varphi(q)N^{1-\frac1n}\lnc^{A+2}}.
\end{align*}
Substituting the expression for the function $\F(\alpha)$, i.e., the right-hand side of the last formula, into (\ref{I(a,q)=intF(alpha;N,H)}), and then using the estimate $S(a,q)\ll q^{1-\frac1n}$ (Lemma \ref{Lemma otsenka polnoy trigsumm}) and relation (\ref{formula poryadok H3}), we find:
\begin{align}
I&(a,q)=\frac{H_3}{2H\ln(\mu_1N)\ln(\mu_2N)}\cdot\frac{\mu^2(q)S(a,q)}{q\varphi^2(q)}e\left(-\frac{aN}q\right)J(H) +\R_5,\label{I(a,q)= ...J(H)+R_5}\\
&J(H)=\int\limits_{|\lambda|\le\lnc^2H^{-1}}\frac{\sin^32\pi\lambda H}{\pi^3\lambda^3}d\lambda,\qquad \R_5\ll\frac{H^2}{q^\frac1n\varphi(q)N^{1-\frac1n}\lnc^{A}}.\nonumber
\end{align}
By replacing $J(H)$ with the corresponding improper integral, which is close to $J(H)$ and independent of $\lnc $, and using the formula (see \cite{Uitiker+Vatson} p.174 )
\begin{align*}
\int\limits_0^\infty\frac{\sin^nmu}{u^n}du=&\frac{\pi
m^{m-1}}{2^n(n-1)!}\left[n^{n-1}-\frac{n}{1!}(n-2)^{n-1}+\frac{n(n-1)}{2!}(n-4)^{n-1}
+\ldots\right],
\end{align*}
 for $m=1$ and $n=3$, we find
 \begin{align*}
J(H)&=\frac{4H^2}{\pi}\int\limits_{|u|\le2\pi\lnc^2}\frac{\sin^3u}{u^3}d\lambda=\frac{8H^2}{\pi}\int\limits_0^\infty\frac{\sin^3u}{u^3}du +O\left(\frac{H^2}{\lnc^6}\right)=3H^2+O\left(\frac{H^2}{\lnc^6}\right).
\end{align*}
Substituting the value of the integral $J(H)$ into formula (\ref{I(a,q)= ...J(H)+R_5}), we find
\begin{align*}
&I(a,q)=\frac{3HH_3}{2\ln(\mu_1N)\ln(\mu_2N)}\cdot\frac{\mu^2(q)S(a,q)}{q\varphi^2(q)}e\left(-\frac{aN}q\right)+\R_6,\\
&\R_6\ll\frac{\mu^2(q)|S(a,q)|}{q\varphi^2(q)}\cdot\frac{HH_3}{\lnc^8}+\R_5\ll\frac{H^2}{q^\frac1n\varphi^2(q)N^{1-\frac1n}\lnc^8} +\frac{H^2}{q^\frac1n\varphi(q)N^{1-\frac1n}\lnc^A}.
\end{align*}
Substituting the obtained value of the integral $I(a,q)$ into (\ref{predst I(M_1)= sum sumI(a,q)}), and then, when estimating $\R_7$, using the condition $A=3+\frac{n-1}n\eta(n)$, we get
\begin{align}
I(\M_1)&=\frac{3HH_3}{2\ln(\mu_1N)\ln(\mu_2N)}\sum_{q\le\lnc^{\eta(n)}}\frac{\mu^2(q)}{q\varphi^2(q)}\sum_{\substack{a=0\\ (a,q)=1}}^qS(a,q)e\left(-\frac{aN}{q}\right)+\R_7.\label{predst I(M_1)= sum sumI(a,q)+O(...)}\\
\R_7&\ll\frac{H^2}{N^{1-\frac1n}}\sum_{q\le\lnc^{\eta(n)}}\sum_{\substack{a=0\\ (a,q)=1}}^{q-1}\left(\frac1{q^\frac1n\varphi^2(q)\lnc^8} +\frac1{q^\frac1n\varphi(q)\lnc^A}\right)\ll \nonumber\\
&
\ll\frac{H^2}{N^{1-\frac1n}\lnc^8}+\frac{H^2}{N^{1-\frac1n}\lnc^{A-\frac{n-1}n\eta(n)}}\ll\frac{H^2}{N^{1-\frac1n}\lnc^3}.\nonumber
\end{align}
The sum over $q$ in (\ref{predst I(M_1)= sum sumI(a,q)+O(...)})  can be replaced by a closely related infinite series, independent of the power of $\lnc $:
\begin{align}\label{pradst ...=Si (N)-R(N)}
&\sum_{q\le \lnc^{\eta(n)}}\frac{\mu^2(q)}{q\varphi^2(q)}\Phi(q,N)=\Si (N)-R(N),\qquad \Phi(q,N)=\sum_{\substack{a=0\\ (a,q)=1}}^qS(a,q)e\left(-\frac{aN}{q}\right),\\
&\Si (N)=\sum_{q=1}^{\infty }\frac{\mu^2(q)}{q\varphi^2(q)}\Phi(q,N),\qquad R(N)=\sum_{q>\lnc^{\eta(n)}}\frac{\mu^2(q)}{q\varphi^2(q)}\Phi(q,N). \nonumber
\end{align}
We represent the special series $\Si(N)$ as an infinite product over all prime numbers. To do so, we first show that the sum $\Phi(q,N)$ is a multiplicative function. Let $q=q_1q_2$, with $(q_1,q_2)=1$. Then, representing the summation variable $a$ as
$$
a=a_1q_2+a_2q_1, \quad (a_1,q_1)=1, \quad 1\le a_1\le q_1,\quad (a_2,q_2)=1, \quad 1\le a_2\le q_2,
$$
we find
\begin{align}\label{predst Phi(q_1q_2)=sumsumS(a_1q_2+a_2q_1,q_1q_2)e(...)}
&\Phi(q_1q_2,N)=\sum_{\substack{a_1=1\\ (a_1,q_1)=1}}^{q_1} \sum_{\substack{a_2=1\\ (a_2,q)=1}}^{q_2}S(a_1q_2+a_2q_1,q_1q_2)e\left(-\frac{(a_1q_2+a_2q_1)N}{q_1q_2}\right).
\end{align}
Representing the summation variable $x$ in the sum $S(a_1q_2+a_2q_1,q_1q_2)$ as
$$
x=x_1q_2+x_2q_1, \qquad 1\leq x_1\leq q_1,  \quad 1\leq x_2\le q_2,
$$
and noting that
\begin{align*}
(x_1q_2+x_2q_1,q_1q_2)=(x_1q_2+x_2q_1,q_1)(x_1q_2+x_2q_1,q_2)=(x_1q_2,q_1)(x_2q_1,q_2)=(x_1,q_1)(x_2,q_2),
\end{align*}
we obtain
\begin{align*}
&S(a_1q_2+a_2q_1,q_1q_2)=\sum_{x=1}^{q_1q_2}e\left(\frac{(a_1q_2+a_2q_1)x^3}{q_1q_2}\right)= \sum_{x_1=1}^{q_1}\sum_{x_2=1}^{q_2}e\left(\frac{(a_1q_2+a_2q_1)(x_1q_2+x_2q_1)^3}{q_1q_2}\right)=\\
&=\sum_{x_1=1}^{q_1}e\left(\frac{a_1(x_1q_2)^3}{q_1}\right)\sum_{x_2=1}^{q_2}e\left(\frac{a_2(x_2q_1)^3}{q_2}\right) =\sum_{x_1=1}^{q_1}e\left(\frac{a_1x_1^3}{q_1}\right)\sum_{x_2=1}^{q_2}e\left(\frac{a_2x_2^3}{q_2}\right)=S(a_1,q_1)S(a_2,q_2).
\end{align*}
Substituting this equality into the right-hand side of (\ref{predst Phi(q_1q_2)=sumsumS(a_1q_2+a_2q_1,q_1q_2)e(...)}), we obtain
\begin{align*}
\Phi(q_1q_2,N)=&\sum_{\substack{a_1=1\\ (a_1,q_1)=1}}^{q_1}S(a_1,q_1)e\left(-\frac{a_1N}{q_1}\right)\sum_{\substack{a_2=1\\ (a_2,q)=1}}^{q_2}S(a_2,q_2) e\left(-\frac{a_2N}{q_2}\right)=\Phi(q_1,N)\Phi(q_2,N).
\end{align*}
Using the absolute convergence of  $\Si (N)$ and the multiplicativity of $\Phi(q,N)$, we find
$$
\Si(N)=\sum_{q=1}^{\infty}\frac{\mu^2(q)}{q\varphi^2(q)}\Phi(q,N)=\prod_p \left(1+\frac{\Phi(p,N)}{p(p-1)^2}\right).
$$
Now, we calculate $\Phi(p,N)$:
\begin{align}
\Phi(p,N)&=\sum_{a=1}^{p-1}S(a,p)e\left(-\frac{aN}{p}\right)=\sum_{x=1}^p\sum_{a=1}^{p-1}e\left(\frac{a(x^n-N)}{p}\right)=\nonumber\\
&=\sum_{x=1}^p\sum^p_{a=1}e\left(\frac{a(x^n-N)}{p}\right)-p=p(\rho(N,p)-1).\label{formula Phi(p,N)-toch znach}
\end{align}
where  $\rho (N,p)$ is the number of solutions to the congruence $x^n\equiv N\pmod p$.
Thus,
$$
\Si(N)=\prod_p \left(1+\frac{\rho(N,p)-1}{(p-1)^2}\right).
$$

Now, we estimate $R(N)$. Using the fact that $\Phi(q,N)$ is a multiplicative function, $q$ is square-free, formula (\ref{formula Phi(p,N)-toch znach}) holds,  and the estimate $\rho(N,p)\le n$, we obtain
\begin{align*}
R(N)&=\sum_{q>\lnc^{\eta(n)}}\frac{\mu^2(q)}{q\varphi^2(q)}\prod_{p\backslash q}\Phi(p,N)=
\sum_{q>\lnc^{\eta(n)}}\frac{\mu^2(q)}{\varphi^2(q)}\prod_{p\backslash q}(\rho(N,p)-1)\le\\
&\le\sum_{q>\lnc^{\eta(n)}}\frac{\mu^2(q)}{\varphi^2(q)}(n-1)^{\omega(q)}=
\sum_{q>\lnc^{\eta(n)}}\frac{\mu^2(q)}{\varphi^2(q)}\exp\left(\omega(q)\ln(n-1)\right),
\end{align*}
where $\omega (q)$ is the number of distinct prime divisors of $q$. Using the well-known inequalities
$$
\frac{\varphi (q)}{q}\ge \frac{c_\varphi}{\ln\ln q},\qquad \omega(q)\le\frac{c_\omega\ln q}{\ln\ln q},
$$
and then applying the definition of the parameter  $\eta(n)$, we find
\begin{align*}
R(N)&\ll\sum_{q>\lnc^{\eta(n)}}\frac{(\ln\ln q)^2}{c_\varphi^2q^2}\exp\left(\frac{c_\omega\ln(n-1)\ln q}{\ln\ln q}\right)\ll\frac1{\lnc^{0,9\eta(n)}}\ll\frac1{\lnc^2}.
\end{align*}
Thus, the relation (\ref{pradst ...=Si (N)-R(N)})  takes the form
\begin{align*}
\sum_{q\le \lnc  ^{736}}\frac{\mu^2(q)}{q\varphi^2(q)}\Phi(q,N)&=\Si (N)+O\left(\lnc^{-2}\right).
\end{align*}
Substituting the right-hand side of this equality into (\ref{predst I(M_1)= sum sumI(a,q)+O(...)}), and then applying formula (\ref{formula poryadok H3})  and the relation
\begin{align*}
\frac1{\ln(\mu_kN)}-\frac1{\lnc}&=\frac{-\ln\mu_k}{(\lnc-\ln\mu_k)\lnc}\ll\frac1{\lnc^2},
\end{align*}
 we obtain
$$
I(\M_1)=\frac{3\Si(N)H^2}{n\mu_3^{1-\frac1n}N^{1-\frac1n}\lnc^2}+O\left(\frac{H^2}{N^{1-\frac1n}\lnc^3}\right).
$$

\subsection{Estimation of the integral $I(\M_2)$}
We have
$$
I(\M_2)=\int_{\M_2}e(-\alpha N)\S(\alpha;N_1,2H)\S(\alpha;N_2,2H)T(\alpha;N_3,H_3)d\alpha
$$
Proceeding to the estimates and applying the Cauchy inequality for integrals, we find
\begin{align*}
&I(\M_2)\ll\max_{\alpha \in\M_2}|T(\alpha;N_3,H_3)|\int_0^1|\S(\alpha;N_1,2H)||\S(\alpha;N_2,2H)|d\alpha =\\
&=\max_{\alpha\in\M_2}|T(\alpha;N_3,H_3)|\left(\int_0^1|\S(\alpha;N_1,2H)|^2d\alpha\right)^\frac12 \left(\int_0^1|\S(\alpha;N_2,2H)|^2d\alpha\right)^\frac12 =\\
&=\max_{\alpha\in\M_2}|T(\alpha;N_3,H_3)|\left(\pi\left(\mu_1N+H\right)-\pi\left(\mu_1N-H\right)\right)^\frac12 \left(\pi\left(\mu_2N+H\right)-\pi\left(\mu_2N-H\right)\right)^\frac12.
\end{align*}
Applying to the last two factors on the right-hand side of the obtained formula, the relation
$$
H=H=N^{1-\theta(n)}\lnc^{\omega(n)}\ge \left(\mu_kN+H\right)^{\frac7{12}+\varepsilon},\qquad k=1;\ 2,
$$
and Lemma \ref{Lemma Huxley pi(x)-pi(x-y)}, we find
$$
\pi(\mu_kN+H)-\pi(\mu_kN-H)\ll\frac{H}{\lnc }.
$$
Hence,
\begin{equation}\label{Integral I(M2)}
I(\M_2) \ll\frac{H}\lnc\max_{\alpha \in \M_2}|T(\alpha;N_3,H_3)|= \frac{H^2}{N^{1-\frac1n}\lnc^3}\cdot\frac{N^{1-\frac1n}\lnc^2}H\max_{\alpha \in \M_2}|T(\alpha;N_3,H_3)|.
\end{equation}
Let us estimate $T(\alpha;N_3,H_3)$ for $\alpha$  from the set $\mathfrak{M}_2$. If $\alpha\in\mathfrak{M}_2$, then
$$
\alpha=\frac aq+\lambda,\qquad(a,q)=1,\qquad\frac{\lnc^2}H<|\lambda |\le\frac1{q\tau},\qquad 1\le q\le\lnc^{\eta(n)}.
$$
We consider two possible cases:   $\frac{\lnc^2}H<|\lambda |\le\frac1{2nqN_3^{n-1}}$   and   $\frac1{2nqN_3^{n-1}}<|\lambda |\le\frac1{q\tau}$.

{\bf Case 1.} According to Corollary \ref{Sledst1 lemmi ob poved kor trig summi Weyl} of Lemma \ref{Lemma ob poved kor trig summi Weyl}, we have
\begin{align}\label{formula T3(lambda;N3,H3)}
T(\alpha;N_3,H_3)=\frac{H_3S(a,q)}{q}\gamma(\lambda;N_3,H_3)+O\left(q^{\frac12+\varepsilon}\right).
\end{align}
We now estimate the exponential integral  $\gamma(\lambda;N_3,H_3)$  based on the magnitude of the first derivative (by Lemma  \ref{Lemma otsenka trig integ po pervoy proizv}). Assuming $f(u)=\lambda\left(N_3-\frac{H_3}2+H_3t\right)^n$, the second derivative\linebreak  $f''(u)=n(n-1)\lambda H_1^2\left(N_3-\frac{H_3}2+H_3t\right)^{n-2}$ does not change the sign. Hence the first derivative $f'(u)$ is monotonic and satisfies the inequality
\begin{align*}
|f'(u)|=n|\lambda|H_3N_3^{n-1}\left(1+\frac{H_3}{N_3}\left(t-\frac12\right)\right)^{n-1}\ge n|\lambda|H_3N_3^{n-1}\left(1-\frac{(n-1)H_3}{N_3}\right).
\end{align*}
Applying the following formulas
$$
H_3N_3^{n-1}=2H\left(1+O\left(\frac HN\right)\right),\qquad \frac{H_3}{N_3}=\frac{2H}{n\mu_3N}\left(1+O\left(\frac HN\right)\right),
$$
to the right-hand side of this inequality, which in turn follow from formulas (\ref{formula poryadok N3})  and (\ref{formula poryadok H3}), we obtain
\begin{align*}
|f'(u)|\ge 2n|\lambda|H\left(1+O\left(\frac HN\right)\right) \left(1-\frac{2(n-1)H}{n\mu_3N}+O\left(\frac{H^2}{N^2}\right)\right)\ge n|\lambda|H.
\end{align*}
Thus, by applying Lemma \ref{Lemma otsenka trig integ po pervoy proizv},  for $m=n|\lambda|H$  and $M=1$, and using the condition $|\lambda|>\lnc^2H^{-1}$, we obtain
$$
|\gamma (\lambda;N_3,H_3)|\le\frac1{n|\lambda|H}\le\frac1{n\lnc^2}.
$$
Proceeding to the estimates in (\ref{formula T3(lambda;N3,H3)}), and then using relation (\ref{formula poryadok H3}), the trivial bound $S(a,q)\le q$, the previously obtained estimate for $\gamma(\lambda;N_3,H_3)$, the condition $q\le \lnc^{\eta(n)}$, and formula (\ref{formula opred teta i omega}), we successively obtain
\begin{align}
|T(\alpha;N_3,H_3)|&\ll\frac{H}{N^{1-\frac1n}\lnc^2}+\lnc^{\left(\frac12+\varepsilon\right)\eta(n)}
=\frac{H}{N^{1-\frac1n}\lnc^2}\left(1+\frac{N^{1-\frac1n}\lnc^{\left(\frac12+\varepsilon\right)\eta(n)+2}}H\right)=\nonumber\\ &=\frac{H}{N^{1-\frac1n}\lnc^2}\left(1+N^{-\frac{n-2}{n(n-1)}}\lnc^{\left(\frac12+\varepsilon\right)\eta(n)+2-\omega(n)}\right) \ll\frac{H}{N^{1-\frac1n}\lnc^2}.
\label{M2}
\end{align}

{\bf Case 2.} In this case, when
$$
\frac1{2nqN_3^{n-1}}<|\lambda |\le\frac1{q\tau},
$$
using Corollary \ref{Sledst2 lemmi ob poved kor trig summi Weyl} of Lemma \ref{Lemma ob poved kor trig summi Weyl}, with
$x=N_3$, $y=H_3$, and the relation (\ref{formula poryadok N3}), followed by the explicit value of the parameter $H$, i.e., formula (\ref{formula opred teta i omega}), we have
\begin{align*}
T(\alpha,&N_3,H_3)\ll q^{1-\frac1n}\ln q +\min\left(H_3q^{-\frac1n}, N_3^\frac12q^{\frac12-\frac1n}\right)\ll q\lnc+N_3^\frac12q^{\frac12-\frac1n}\ll\\
&\ll\lnc^{\eta(n)+1}+N^\frac1{2n}\lnc^{\left(\frac12-\frac1n\right)\eta(n)}= \frac{H}{N^{1-\frac1n}\lnc^2}\left(\frac{N^{1-\frac1n}\lnc^{\eta(n)+3}}H+ \frac{N^{1-\frac1{2n}}\lnc^{\left(\frac12-\frac1n\right)\eta(n)+2}}{H}\right)=\\
&=\frac{H}{N^{1-\frac1n}\lnc^2}\left(N^{-\frac{n-2}{n(n-1)}}\lnc^{\eta(n)+3}+ N^{-\frac{n-3}{2n(n-1)}}\lnc^{\varkappa(n)}\right),\quad \varkappa(n)=\left(\frac12-\frac1n\right)\eta(n)+2-\omega(n).
\end{align*}
From this, noting that $\varkappa(3)=-\frac{14}3$, we obtain
$$
T(\alpha,N_3,H_3)\ll\frac{H}{N^{1-\frac1n}\lnc^2}.
$$
From this estimate and from estimate (\ref{M2}), in view of the relation (\ref{Integral I(M2)}) , for all $\alpha\in \M_2$ we get
\begin{align*}
I(\M_2) \ll\frac{H^2}{N^{1-\frac1n}\lnc^3}.
\end{align*}

\subsection{Estimation of the integral $J(\mathfrak{m})$}
Proceeding similarly to the estimation of $J(\mathfrak{M}_2)$, we obtain
\begin{align}
J(\mathfrak{m})\ll\frac{H^2}{N^{1-\frac1n}\lnc^3}\cdot\frac{N^{1-\frac1n}\lnc^2}H\max_{\alpha\in\mathfrak{m}}|T(\alpha;N_3,H_3)|. \label{formula I(m)-1}
\end{align}
Now, let us estimate $T(\alpha;N_3,H_3)$ for $\alpha$  from the set $\mathfrak{m}$. If $\alpha\in\mathfrak{m}$, then
$$
\alpha=\frac aq+\lambda,\qquad (a,q)=1, \qquad |\lambda|\le\frac1{q\tau},\qquad \lnc^{\eta(n)}<q\le\tau,\qquad \eta(n)=2^{n+1}+(n-1)^2.
$$
Using Theorem \ref{Teor otsenka kor summ Weyl n>=3}, with $x=N_3$, $y=H_3$, and using the relations
$$
H_3\asymp\frac H{N^{1-\frac1n}},\qquad N_3\asymp N^\frac1n,\qquad  \tau\frac H{N^\frac1n},
$$
which follow from formulas (\ref{formula poryadok H3}), (\ref{formula poryadok N3})  and (\ref{formula poryadok tau}), and then using the explicit value of the parameter $H$, i.e., formula (\ref{formula opred teta i omega}), we have
\begin{align*}
T(\alpha;&N_3,H_3)|\le 2H_3\left(\frac1q+\frac1{H_3}+\frac{q}{H_3^n}\right)^{2^{-n}}(\ln qH_3)^{(n-1)^22^{-n}}\ll\\
&\ll\frac H{N^{1-\frac1n}}\left(\frac1{\lnc^{2^{n+1}+(n-1)^2}}+\frac{N^{1-\frac1n}}H+ \frac{N^{n-1-\frac1n}}{H^{n-1}}\right)^{2^{-n}}\lnc^{(n-1)^22^{-n}}=\\
&=\frac H{N^{1-\frac1n}\lnc^2}\left(1+\frac{N^{1-\frac1n}\lnc^{2^{n+1}+(n-1)^2}}H+ \left(\frac{N^{1-\frac1{(n-1)n}}\lnc^{\frac{2^{n+1}}{n-1}+n-1}}H\right)^{n-1}\right)^{2^{-n}}=\\
&=\frac H{N^{1-\frac1n}\lnc^2}\left(2+N^{-\frac{n-2}{n(n-1)}}\lnc^{\frac{(n-2)2^{n+1}}{n-1}+n^2-3n+2}\right)^{2^{-n}}
\ll\frac H{N^{1-\frac1n}\lnc^2}.
\end{align*}
Substituting this estimate into (\ref{formula I(m)-1}), we obtain
$$
J(\mathfrak{m})\ll\frac{H^2}{N^{1-\frac1n}\lnc^3}.
$$
The theorem is proven.

\end{document}